\documentclass[a4paper]{article}

\usepackage{amsmath,amsthm,amssymb}
\usepackage{hyperref}
\usepackage{microtype}
\usepackage{xcolor}
\usepackage{algorithm}
\usepackage{algpseudocode}

\hypersetup{
	colorlinks=true,
	linkcolor=black,
	citecolor=blue,
	urlcolor=blue}

\newtheorem{thm}{Theorem}[section]
\newtheorem{lma}[thm]{Lemma}
\newtheorem{cor}[thm]{Corollary}

\theoremstyle{definition}
\newtheorem{dfn}[thm]{Definition}

\newcommand{\R}{\mathbb{R}}
\newcommand{\N}{\mathbb{N}}

\DeclareMathOperator{\supp}{supp}
\newcommand{\calP}{\mathcal{P}}

\begin{document}

\title{Convergence proof for the GenCol algorithm \\
	in the case of two-marginal optimal transport}

\author{Gero Friesecke \and Maximilian Penka\thanks{MP was partially funded by the Deutsche Forschungsgemeinschaft (DFG, German Research Foundation) through the Collaborative Research Center TRR 109 “Discretization in Geometry and Dynamics“, Projektnummer 195170736}}

\date{\it Department of Mathematics, TUM}

\maketitle

\begin{abstract}
The recently introduced Genetic Column Generation (GenCol) algorithm has been numerically observed to efficiently and accurately compute high-dimensional optimal transport plans for general multi-marginal problems, but theoretical results on the algorithm have hitherto been lacking. The algorithm solves the OT linear program on a dynamically updated low-dimensional submanifold consisting of sparse plans. The submanifold dimension exceeds the sparse support of optimal plans only by a fixed factor $\beta$. Here we prove that for $\beta \geq 2$ and in the two-marginal case, GenCol always converges to an exact solution, for arbitrary costs and marginals. The proof relies on the concept of c-cyclical monotonicity. As an offshoot, GenCol rigorously reduces the data complexity of numerically solving two-marginal OT problems from $O(\ell^2)$ to $O(\ell)$ without any loss in accuracy, where $\ell$ is the number of discretization points for a single marginal.
At the end of the paper we also present some insights into the convergence behavior in the multi-marginal case.
\end{abstract}

\section{Introduction}\label{S:Intro}
Large-scale discrete optimal transport problems are difficult to solve numerically because the full problem has a huge number of possible configurations. At the same time it is guaranteed that a rather sparse solution exists, a particularly well known manifestation in continuous OT being Brenier's theorem. This effect is especially important in the multi-marginal case, but occurs already in the classical two-marginal case when the support size of both marginals is large.

In recent years, computational strategies for optimal transport were driven by the idea of approximating the problem by adding an entropy-like penalty term. This transforms the problem into a strictly convex and more robust optimization problem, which can be solved in short time  using the Sinkhorn algorithm as long as the overall number of unknowns remains moderate. In theory this approach, called entropic optimal transport (EOT), is also valid for two-marginal problems in high dimension or general multi-marginal problems.

Unfortunately, this approach corresponds to smearing out the transport plan, yielding a huge amount of configurations in its support: The true optimizer of the EOT problem has the same support as the full product measure of its marginals. Thus the support size scales polynomially in the support size of the marginals, and exponentially in the number of marginals. Recent approaches alleviate this problem by truncation or multi-scale methods \cite{schmitzer2019stabilized} or -- most recently -- low-rank approximation \cite{scetbon2021low,strossner2022low}.

An alternative approach was proposed in \cite{friesecke2022genetic, friesecke_penka2022}. Rather than looking for further refinements of interior point methods, which struggle to solve huge programs, one goes back to the original linear program, and exploits that the OT program possesses extremely sparse solutions. If the $i$-th marginal is supported on $\ell_i$ points, there exist optimal plans with support size less than the {\it sum} of the $\ell_i$, rather than their product \cite{friesecke_penka2022}. 

The standard approach to solve linear programs is the simplex algorithm and its descendants, most promisingly {\it Column Generation} \cite{lubbecke2005selected}. The latter tackles large LPs by iteratively solving smaller (``reduced'') problems on a dynamically evolving subset of all variables. In optimal transport, every variable of the LP corresponds to a possible configuration in the product of the marginal domains, so Column Generation corresponds to solving the OT problem on a subset of the product of the marginal domains. Thus Column Generation can in principle exploit sparsity and  find the exact optimal plan of the full problem, {\it provided} a sparse superset of its support is known. But in practice such a superset has to be found.

There are two obstructions for doing so. First, in Column Generation the generation of new variables is originally done by solving a second optimization problem, the so called pricing problem. Unfortunately the pricing problem for large problems is still expensive; in the multi-marginal case it has been proven to be NP hard \cite{altschuler2021, friesecke2022genetic}. For LPs of moderate size this problem can be alleviated by generating new configurations at random, but in this randomized version one needs to try many configurations, again leading to an unacceptable slowdown for large problems. 
Second, the algorithm increases the size of the LP in each iteration step, making the iterations slower and slower and lacking any convergence guarantee until the size has reached the original LP size that one wanted to avoid! 

For these reasons, \cite{friesecke2022genetic} invented 
\begin{itemize}
    \item a genetic search rule, restricting the number of possible proposals from all configurations to an update of one entry in one active configuration
    \item a genetic tail-clearing rule which discards those configurations which have been inactive the longest, to keep the overall support size at a fixed small multiple of the size of sparse optimizers.
\end{itemize}

The resulting algorithm, which was termed Genetic Column Generation (GenCol), exhibited in several application examples of interest a spectacularly accelerated convergence to global optimizers. A theoretical explanation has hitherto been lacking.

In this paper we present a rigorous proof of convergence to a global optimizer in the case of two marginals.
The fact that GenCol cannot get stuck in a local optimizer is far from obvious since the upper bound on the support size makes the reduced problem non-convex.
The proof relies on the concept of $c$-cyclical  monotonicity which is well known in the theory of optimal transport. It finds here a beautiful application and yields an intuitive understanding of the algorithm.

Our arguments, while rigorously and non-trivially reducing the storage cost, do not yield  a rigorous bound on the convergence speed,
which -- in numerical examples  -- is observed to be exponential \cite{friesecke2022genetic, friesecke_penka2022}.
By contrast, for Sinkhorn as well as some classical LP algorithms requiring  access to the full state space, the convergence speed has been rigorously estimated.
See \cite{franklin1989scaling} for Sinkhorn,
\cite{li2022maxflow} for max-flow min-cut,
and \cite{chandrasekaran2016cutting} for cutting plane applied to perfect matching. Let us also mention the numerical study \cite{dong2020study} which compares the runtimes of some standard two-marginal OT algorithms.

In section \ref{S:MMOT} we analyze the multi-marginal case, for which the GenCol algorithm was originally proposed. We rigorously justify convergence of the algorithm to a global optimizer {\it provided} the search rule finds any possible configuration with positive probability. Hence with such a search rule, GenCol rigorously reduces the storage complexity from exponential to linear in the number of variables. However the price to pay is that the number of search steps might be exponentially large.
By contrast, GenCol with the efficient search rule from \cite{friesecke2022genetic, friesecke_penka2022} - with its one-entry-at-a-time update which requires only quadratically many search steps in the number of marginals - might fail to converge to a global optimizer, at least for general costs. See section \ref{S:Counterexample} for a counterexample.
The design and analysis of updating rules for the multi-marginal case which are both efficient and yield rigorous global convergence for practically relevant costs 
is an interesting open question.

\section{c-Cyclical Monotonicity}\label{S:cCM}
Given two probability measures $\mu_1,\mu_2$ on Polish spaces $X$ respectively $Y$, the optimal transport problem is the following:
\begin{align*}
    \begin{split}
    \operatorname{minimize}&\quad \mathcal{F}[\gamma] := \int_{X\times Y} c(x,y)\,d\gamma(x,y)\quad \text{over }\gamma\in \calP(X\times Y)\\
        \text{subject to}&\quad \begin{cases}
        \gamma(A\times Y) = \mu_1(A) & \text{for all measurable } A \subset X\\
        \gamma(X\times B) = \mu_2(B) & \text{for all measurable } B \subset Y.
        \end{cases} 
    \end{split}  
\end{align*}
where $\calP$ denotes the set of probability measures. Solutions to the constraints are called transport plans. Optimality of a transport plan $\gamma$ can be characterized by a condition on its support, called $c$-cyclical monotonicity.

\begin{dfn}[see e.g. \cite{Santambrogio2015}, Def. 1.36]
Given a function $c\colon X\times Y \to \R \cup \{+\infty\}$, we say that a set $\Gamma \subset X \times Y$  is \emph{c-cyclically monotone} (c-CM) if for every $k\in \N$, every permutation $\sigma: \{1,\dots,k\} \to \{1,\dots,k\}$, and every finite set of points $\{(x_1,y_1),...,(x_k,y_k)\} \subset \Gamma$ we have
\begin{equation*}
    \sum_{i=1}^k c(x_i,y_i) \leq \sum_{i=1}^k c(x_i,y_{\sigma(i)}).
\end{equation*}
\end{dfn}
While it is easy to see (at least in the discrete case) that this is a necessary condition on the support of an optimal plan, it turns out to also be sufficient.

\begin{thm}[see e.g. \cite{Santambrogio2015}, Thm. 1.49] \label{thm:c-CM}
Suppose $X$ and $Y$ are Polish spaces and $c\colon X\times Y \to \R$ is uniformly continuous and bounded. Given $\gamma \in \mathcal{P}(X\times Y)$, if $\operatorname{supp}(\gamma)$ is c-CM then $\gamma$ is an optimal transport plan between its marginals $\mu_1 = (\pi_1)_\sharp\gamma$ and $\mu_2 = (\pi_2)_\sharp\gamma$ for the cost $c$.
\end{thm}

\section{Sparsity of optimal plans}\label{S:sparsity}
The support of optimal transport plans is typically a much smaller set than the product of the supports of the marginals. Rather than going into classical variants for convex costs like Brenier's theorem and their interesting relation to c-cyclical monotonicity, we focus directly on a discrete version for general costs which informed the design of the GenCol algorithm and is useful for its analysis. 

For $X$ and $Y$ discrete, $|X| = \ell_1, |Y| = \ell_2$, the objective function $\mathcal F$ becomes a finite sum and the OT problem a linear program in standard form:
\begin{align}\label{OT} \tag{OT}
    \begin{split}
        \operatorname{minimize}\quad &\langle c,\gamma\rangle := \sum_{(x,y) \in X\times Y} c(x,y)\gamma(x,y) \text{ over } \gamma: X\times Y \to [0,\infty)\\
        \text{subject to}\quad &\gamma \in \Pi(\mu_1,\mu_2) :\Leftrightarrow \begin{cases} \sum_{y \in Y} \gamma(x_0,y) = \mu_1(x_0)\,\forall x_0 \in X \\ \sum_{x \in X} \gamma(x,y_0) = \mu_2(y_0)\, \forall y_0 \in Y,\end{cases}
    \end{split}
\end{align}
where the measures $\mu_1,\mu_2$ and $\gamma$ were identified with their densities with respect to the counting measures on their domains.

\begin{thm} \label{T:sparse}
  Suppose $X$ and $Y$ are discrete with $|X|=\ell_1$, $|Y|=\ell_2$. Then any extreme point of the Kantorovich polytope $\Pi(\mu_1,\mu_2)$ is supported on at most $\ell_1+\ell_2-1$ points. In particular, for any cost $c\, : \, X\times Y\to \R$ and any marginals, the OT problem \eqref{OT} possesses an optimizer supported on at most $\ell_1+\ell_2-1$ points. 
\end{thm}

This can be deduced from well known results on extremal solutions in linear programming. For a self-contained and simple proof using geometry of convex polytopes see \cite{friesecke_penka2022}.

\section{Genetic Column Generation Algorithm} \label{S:GenCol}
The algorithm doesn't deal with the full OT problem but only its restrictions to certain subsets of $X\times Y$ whose size is of the order of the support size of optimizers from Theorem \ref{T:sparse}.  

We call a subset $\Omega \subset  X\times Y$ a {\it feasible subset of configurations} if $\Pi(\mu,\nu) \cap \{\gamma : \supp(\gamma) \subset \Omega\}$ is non-empty. Given such a subset, we define the \emph{reduced problem} to be
\begin{align}\label{ROT} \tag{ROT}
    \begin{split}
        \operatorname{minimize}\quad &\langle c,\gamma\rangle \text{ over } \gamma: X\times Y \to [0,\infty)\\
        \text{subject to}\quad &\gamma \in \Pi(\mu_1,\mu_2) \\
        &\text{ and } \supp(\gamma)\subseteq\Omega.
    \end{split}
\end{align}
Because $\gamma$ is a discrete measure, $\supp(\gamma)$ is the set of all $(x,y)\in X\times Y$ with $\gamma(x,y)\neq 0$. Thus the reduced problem amounts to reducing the variables in the linear program to the values of $\gamma$ on configurations in $\Omega$ (and setting the values outside $\Omega$ to zero), and not changing the constraints. As the values outside $\Omega$ no longer need to be considered, this shrinks the size of the program to that of $\Omega$.

Before we come to genetic column generation, let us describe  {\it classical column generation}. Unlike genetic column generation it does not restrict the size of $\Omega$, and works as follows. Given a feasible initial set $\Omega$, the first step is to solve the reduced problem. The second step is to generate a new configuration $(x',y') \notin \Omega$ which is added to $\Omega$ and improves the solution. The two steps are iterated until no more improving configurations exist.

The second step relies on the dual of the reduced problem \eqref{ROT}, 
\begin{align}
    \tag{D-ROT} \label{D-ROT}
    \begin{split}
        \operatorname*{maximize}\quad &\langle \mu_1,u_1\rangle + \langle \mu_2,u_2\rangle \text{ over}\, u_1 : X \to \R, u_2 : Y \to \R \\
        \text{such that}\quad &u_1(x) + u_2(y) \leq c(x,y) \quad  \forall (x,y) \in \Omega.
    \end{split}
\end{align}
In comparison, the dual of the full problem \eqref{OT} has the same objective function, but more constraints:
\begin{align}
    \tag{D-OT} \label{DOT}
    \begin{split}
        \operatorname*{maximize}\quad &\langle \mu_1,u_1\rangle + \langle \mu_2,u_2\rangle \text{ over}\, u_1 : X \to \R, u_2 : Y \to \R \\
        \text{such that}\quad &u_1(x) + u_2(y) \leq c(x,y) \quad  \forall (x,y) \in X\times Y.
    \end{split}
\end{align}
Hence every dual optimizer for the full problem is admissible in the reduced problem \eqref{D-ROT}, but a dual optimizer for the reduced problem might violate a constraint of the full problem \eqref{DOT}. If, however, a dual optimizer for the reduced problem is admissible for  \eqref{DOT} then it is already optimal for \eqref{DOT}:
\begin{lma}
    Let $(\gamma^\star,(u_1^\star,u_2^\star))$ be a pair of optimizers for the reduced problems (\ref{ROT},~\ref{D-ROT}). If $(u_1^\star,u_2^\star)$ is \textbf{admissible} for the dual of the full problem \eqref{DOT}, then $\gamma^\star$ is \textbf{optimal} for \eqref{OT}.
\end{lma}
For a proof of this classical result translated into the present context and language of OT see \cite{friesecke_penka2022}.
Hence new configurations $(x',y')\notin\Omega$ can be sought by checking if they violate the dual constraint of the full problem \eqref{DOT}, i.e. if they satisfy the following acceptance criterion:
\begin{equation}\label{acc}\tag{Acc}
  u_1^\star(x') + u_2^\star(y') - c(x',y') > 0.
\end{equation}
Due to economic interpretations this difference is called gain. In classical column generation this gain is maximized over all configurations, constituting the  so-called {\it pricing problem}. 

The following difficulties arise when applying column generation to large LPs, as already pointed out in the Introduction. (i) The pricing problem is too expensive; and the empirical strategy of instead generating configurations $(x',y')\notin\Omega$ independently at random until one of them satisfies \eqref{acc} requires too many trials, especially in the multi-marginal case. (ii) Regardless of how one searches for new configurations, the subset $\Omega$ 
grows in each iteration step, making the iterations slower and slower and lacking any convergence guarantee until the size
has reached the original LP size that one wanted to avoid.

{\it Genetic column generation} \cite{friesecke2022genetic,friesecke_penka2022} tackles these difficulties as follows.
 
(i) Motivated by machine learning protocols in unsupervised learning, the algorithm first proposes new configurations originating from currently active configurations, i.e. $(x,y) \in \operatorname{supp}(\gamma) \subset \Omega$:
one picks an active configuration at random (``parent''), then proposes an offspring (``child'') by changing one entry of the parent configuration. Explicitly,
\begin{equation}\label{eq:ch2}
    \begin{split}
    &\text{given a parent } (x,y) \in \supp(\gamma), \\[-1mm]
    &\text{pick a random child in } \bigl(\supp(\mu_1)\! \times \!\{y\}\bigr) \;\cup\; \bigl(\{x\}\! \times\! \supp(\mu_2)\bigr).
    \end{split}
\end{equation}
The offspring is then accepted if its gain is positive. The rough analogy to ML is that the proposal step mimics an SGD step and the acceptance mimics learning from an adversary (in this case, the current dual). In fact, the proposal step in the first version of GenCol was even more similar to SGD, in that entries were points on a regular grid and children were proposed from neighbouring sites of parents.

(ii) The size of $\Omega$ is restricted to remain of the order of the support size of optimizers from Theorem \ref{T:sparse}. More precisely, one introduces a  hyperparameter $\beta > 1$ and a tail clearing rule which guarantees that 
\begin{equation}\label{eq:suppbd}
   |\Omega| \le \beta\cdot (\ell_1+\ell_2).
\end{equation}
Tail clearing means that whenever, after accepting a child, $\Omega$ violates \eqref{eq:suppbd}, the oldest unused configurations are removed. In practice, one chooses $3 \lesssim \beta \lesssim 5$ and discards a batch of $\ell_1+\ell_2$ configurations whenever $|\Omega|$ exceeds $\beta\cdot(\ell_1+\ell_2)$. 
The hyperparameter $\beta$ does not depend on the sizes $\ell_1$ and $\ell_2$ of $X$ and $Y$.

See Algorithm \ref{algo:GenCol} for a summary of the algorithm. 

\begin{algorithm} 
\caption{Genetic Column Generation} 
\label{algo:GenCol}
\begin{algorithmic}[1]
\Require Marginals $\mu_1,\mu_2$; feasible set $\Omega$ satisfying \eqref{eq:suppbd}; hyperparameter $\beta > 1$
\While{TRUE}
  \State $(\gamma^\star,u^\star) \gets \text{solution to  \eqref{ROT}, \eqref{D-ROT}}$
  \Repeat
    \State Sample a parent in $\operatorname{supp}(\gamma^\star)$ and a child $(x',y')$ 
  \Until{$u_1^\star(x')+u_2^\star(y') > c(x',y')$ (Acc) {\textbf{or}} all possible offspring were tried}
    \If{$\neg$(Acc)}
      \State \Return $(\gamma^\star,u^\star)$ optimal
    \EndIf
  \State Accept the child: $\Omega \gets \Omega \cup \{(x',y')\}$
  \If{$|\Omega| > \beta \cdot (\ell_1 +\ell_2)$}
    \State remove oldest inactive configurations from $\Omega$
  \EndIf
\EndWhile
\end{algorithmic}
\end{algorithm}

It is not clear why the algorithm should find a global optimum. Can it happen -- due to the tail clearing -- that it instead gets stuck in a local minimum? 

The answer is No, as shown in the next section. Note that (at least in the two-marginal case; see section ... for discussion of the multi-marginal case) every configuration $(x',y') \notin \Omega$ which belongs to the product of the supports of $\mu_1$ and $\mu_2$ is proposed by Algorithm \ref{algo:GenCol} with strictly positive probability, so the genetic proposal of updates is not a restriction. However, the tail-clearing turns the original, convex state space $\Pi(\mu_1,\mu_2)$ into the nonconvex state space $\Pi(\mu_1,\mu_2)\cap\{\gamma\in\calP(X\times Y)\, | \, |\supp \gamma|\le \beta\cdot(\ell_1+\ell_2)$, making the question of global convergence nontrivial.

\section{Convergence}\label{S:convergence}
Before giving the proof of convergence, we must specify line 2 (solving the reduced problem and its dual) and line 4 (sample a parent and a child) of Algorithm \ref{algo:GenCol} more precisely. 

{\it Line 2.} First, in degenerate cases optimal plans may not be unique, so for convergence it is mandatory that $\gamma^*$ in Algorithm \ref{algo:GenCol} is updated only if the previous plan is no longer minimizing. Second, we require the linear programming solver to provide 
a solution $\gamma^*$ which satisfies the support size bound from Theorem \ref{T:sparse}. If one uses the simplex algorithm with a warm start, both these requirements are automatically guaranteed. 

{\it Line 4.} Second, to rigorously implement the second stopping criterion in line 5, one does not sample parents and children in each trial independently, but draws random permutations covering all possibilities  and then tries them one after another until the stopping criterion is satisfied. 

\begin{thm} \label{T:cyclelength} 
  Let $X$ and $Y$ be discrete with $|X|=\ell_1$, $|Y|=\ell_2$, let $c\, : \, X\times Y\to \R$ be any cost, let $\mu_1\in\calP(X)$, $\mu_2\in\calP(Y)$ be any marginals, and let $\Omega\subset X \times Y$ be any feasible subset of configurations. For any optimal solution $\gamma^\star$ for the reduced problem \eqref{ROT} which is an extreme point of the Kantorovich polytope and which is not optimal for the full problem \eqref{OT}, the GenCol proposal and acceptance routine (lines 2--9 of Algorithm \ref{algo:GenCol}) as detailed above finds with positive probability in consecutive steps a superset $\tilde \Omega \supset \Omega$, whose size exceeds that of $\Omega$ by at most $\ell_1+\ell_2-1$ elements, which reduces the total cost:
  \begin{equation*}
    \min_{\gamma \, : \,  \supp(\gamma) \subseteq \tilde \Omega} \mathcal{F}[\gamma] < \min_{\gamma \, : \,  \supp(\gamma) \subseteq \Omega} \mathcal{F}[\gamma].
  \end{equation*}
\end{thm} 

\begin{proof}
Consider an optimal solution $(\gamma^0,u^0)$ for the reduced problem \eqref{ROT} with $\gamma^0$ extremal which is not optimal for the full problem. By Theorem \ref{T:sparse} $\gamma^0$ is sparse with at most $\ell_1+\ell_2-1$ non-zero entries (active configurations). In the following we write $u^0=(u^0_1,u^0_2)$. Due to complementary slackness
\begin{equation*}
  u^0_{{1}}(x) + u^0_{{2}}(y) - c(x,y) = 0 \quad \forall (x,y) \in \supp(\gamma^0). 
\end{equation*}
Because $\gamma^0$ is not optimal for the full problem \eqref{OT}, by Theorem \ref{thm:c-CM} there exists a family $\Gamma = \{(x_1,y_1),...,(x_k,y_k)\}\subset \supp(\gamma^0)$ and a permutation $\sigma \in S_k$ such that
\begin{equation*}
  \sum_{i=1}^k c(x_i,y_i) > \sum_{i=1}^k c(x_i,y_{\sigma(i)}).
\end{equation*}
Note that the bound $|\supp(\gamma)| \leq \ell_1+\ell_2-1$ yields the upper bound $k \leq \ell_1+\ell_2-1$.
Because $\{(x_1,y_1),...,(x_k,y_k)\}\subset \supp(\gamma^0)$, complementary slackness implies
\begin{gather} \addtocounter{equation}{1}
        u^0_1(x_1) + u^0_2(y_1) - c(x_1,y_1) = 0 \tag{\arabic{equation}.1} \label{eq:CS1}\\
         \vdots \notag\\
        u^0_1(x_k) + u^0_2(y_k) - c(x_k,y_k) = 0\tag{\arabic{equation}.k}. \label{eq:CSk}
\end{gather}
After summation,
\begin{align}
    &\;\sum_{i=1}^k \Bigl(u^0_1(x_i) + u^0_2(y_i)\Bigr) - \underbrace{\sum_{i=1}^kc(x_i,y_i)}_{> \sum\limits_{i=1}^k c(x_i,y_{\sigma(i)})} = 0 \label{eq:SumCS}\\
    \Longrightarrow \quad &\;\sum_{i=1}^k\Bigl( u^0_1(x_i) + u^0_2(y_{\sigma(i)}) - c(x_{i},y_{\sigma(i)})\Bigr) > 0 \label{eq:SumCSperm}\\
    \Longrightarrow \quad 
    &\max_{i\in\{1,...,k\}} \left\{u^0_1(x_i) + u^0_2(y_{\sigma(i)}) - c(x_i,y_{\sigma(i)}) \right\} > 0 \label{eq:MaxCS}.
\end{align}
Because, in the two-marginal case, all configurations $(x',y')\notin\Omega$ are proposed by GenCol with positive probability, the element of $\Gamma$ where the maximum in \eqref{eq:MaxCS} is realized, let us call it  $(x_{i_1},y_{\sigma(i_1)})$, is proposed with positive probability, and accepted. 
In the next step the reduced OT problem is resolved on $\Omega$ extended by this element, yielding a new optimal pair $(\gamma^1,(u_1^1,u_2^1))$ and two cases.

Case 1: The optimal plan changes: $\gamma^1 \neq \gamma^0$. In that case, due to the rule that the plan only changes when it must, the optimal cost decreases and we are done. 

Case 2: The optimal plan does not change, $\gamma^1 = \gamma^0$. But the dual solution $(u^1_1,u^1_2)$ must have changed. Because $\gamma^1 = \gamma^0$, we have $(x_1,y_1),\dotsc,(x_k,y_k) \in \supp(\gamma^1$) and eqs. \eqref{eq:CS1} - \eqref{eq:CSk} still hold true with $u^0_1$, $u^0_2$ replaced by $u^1_1$, $u^1_2$. 
But now, since $u^1$ must satisfy the dual constraints on the enlarged configuration set, we also have
\begin{equation*}
u^1_1(x_{i_1}) +  u^1_2(y_{\sigma(i_1)}) - c(x_{i_1},y_{\sigma(i_1)}) \leq 0.
\end{equation*}
Since eqs. \eqref{eq:SumCS}--\eqref{eq:MaxCS} all are also still true with $u^0_1$, $u^0_2$ replaced by $u^1_1$, $u^1_2$, we conclude that 
\begin{equation} \label{eq:maxCS2}
  \max_{i\in \{1,...,k\}\backslash\{i_1\} } \left\{u^1_1(x_i) + u^1_2(y_{\sigma(i)}) - c(x_i,y_{\sigma(i)})\right\} > 0.
\end{equation}
But in the next step, again either the optimal plan changes or the element of $\Gamma$ realizing the maximum in \eqref{eq:maxCS2} is proposed with positive probability and accepted, and so on. After $k$ enlargement steps of $\Omega$, either a change of optimal plan has occurred in some step, or all elements of $\Gamma$ have been accepted with positive probability. But then Case 1 occurs, since the mass $\min\{\gamma^0(x_i,y_i), i = 1,\dotsc,k\} > 0$ can be moved from $\{(x_i,y_i)\}_{i=1}^k$ to $\{(x_i,y_{\sigma(i)}\}_{i=1}^k$, decreasing the total cost.
\end{proof}
One can alternatively see that Case 1 occurs once all elements of $\Gamma$ have been accepted by considering the dual solution $u^k$: otherwise we would have
\begin{equation*}
   u_1^k(x_i) + u_2^k(y_{\sigma(i)}) - c(x_i,y_{\sigma(i)}) \le 0 \; \forall i\in\{1,...,k\}, 
\end{equation*}
but on the other hand \eqref{eq:SumCSperm} must hold with $u^0_1$, $u^0_2$ replaced by $u^k_1$, $u^k_2$, a contradiction.

Convergence of Algorithm \ref{algo:GenCol} now follows as an easy consequence.

\begin{cor} 
  Suppose $X$ and $Y$ are discrete spaces of finite cardinality, and the hyperparameter $\beta$ is $\ge 2$. For any marginals, any cost function, and any feasible initial set $\Omega\subset X\times Y$, GenCol converges with probability 1 to an exact solution of the OT problem \eqref{OT}. 
\end{cor}
\begin{proof}
  The total cost is monotonically decreasing. Moreover $X\times Y$ is finite and every non-optimal plan is improved with positive probability:  Since $\beta\ge 2$, after a tail-clearing the algorithm allows to add more than $\ell_1+\ell_2-1$  configurations, and by Theorem \ref{T:cyclelength} this suffices to  find a plan with lower cost. Therefore the algorithm converges with probability 1.
\end{proof} 

\subsubsection*{Remarks.}
\begin{enumerate}
    \item Note the generality of the cost function.
    \item Shorter families of new configurations are found with higher probability.
    \item Longer tails increase the probability to find also long families, but slow down the simplex algorithm to solve the LP. In practice, a value slightly larger than the minimal value from theory (e.g. $\beta=3$)  works well, and was in fact used on empirical grounds in \cite{friesecke_penka2022}.
\end{enumerate}

\section{The multi-marginal case}\label{S:MMOT}
The algorithm was originally introduced for multi-marginal problems \cite{friesecke2022genetic,friesecke_penka2022}, and its adaptation to this case is straightforward.
One now has $N$ marginals $\mu_1,\dotsc,\mu_N$ on discrete spaces $X_1,\dotsc,X_N$ of sizes $\ell_1,\dotsc, \ell_N$. Plans are nonnegative functions on the product space $X_1\times \dotsb \times X_N$ to $[0,\infty)$, and one seeks to
\begin{align}\label{MMOT} \tag{MMOT}
    \begin{split}
        \operatorname{minimize}\quad &\langle c,\gamma\rangle := \sum_{(x_1,...,x_N) \in X_1\times ... \times X_N} c(x_1,...,x_N)\gamma(x_1,...,x_N)\\
        \text{subject to}\;\; &\gamma \in \Pi(\mu_1,...,\mu_N). 
    \end{split}
\end{align}
This problem possesses an optimizer $\gamma \, : \, X_1\times ... \times X_N \to\R$ supported on at most 
$1+\sum_{i=1}^N (\ell_i - 1)$ points.
Starting from a feasible set of configurations $\Omega \subset X_1\times \dots \times X_N$ satisfying $|\Omega|\le \beta\cdot(\ell_1+\ldots+\ell_N)$, sampling of new configurations works as before: one picks an active configuration (parent) and proposes a child related to the parent,
which is accepted if the gain $\sum_{i=1}^N u_i(x_i)-c(x_1,...,x_N)$ is positive, where $u=(u_1,\ldots,u_N)$, $u_i \, : \, X_i\to\R$, is the current dual solution. Tail clearing is carried out whenever $|\Omega|$ exceeds $\beta\cdot (\ell_1+...+\ell_N)$. Once a child has been accepted and tail clearing has been carried out if necessary, the reduced primal and dual problems on $\Omega$ are re-solved.

Two obvious generalizations of the search rule for children suggest themselves. Either children are proposed by fixing all but one entry or changing all but one entry of the parent configuration:
\begin{gather}
\begin{split}\label{eq:ch_single}
   &\text{Given a parent } (x_1,\dots,x_N) \in \supp(\gamma), \text{ pick a random child in}\\[-1mm]
    &\bigcup\limits_{i=1}^N {\{x_1\} \times \dots \times \{x_{i-1}\} \times \supp(\mu_i) \times \{x_{i+1}\} \times \dots \times \{x_N\}}
\end{split}
\end{gather}
or
\begin{gather}
\begin{split}\label{eq:ch_many}
    &\text{Given a parent } (x_1,\dots,x_N) \in \supp(\gamma), \text{ pick a random child in}\\[-1mm]
    &\bigcup_{i=1}^N {\supp(\mu_1) \times \dots \times \supp(\mu_{i-1}) \times \{x_i\} \times \supp(\mu_{i+1}) \times \dots \times \supp(\mu_N)}.
\end{split}
\end{gather}
For $N=2$ both rules reduce to \eqref{eq:ch2}.

In practice, GenCol with the search rule \eqref{eq:ch_single} turned out to be an extremely fast and accurate method to solve high-dimensional OT problems.
In various test examples with up to $\sim 10^{30}$ variables, it converged to a global optimum of the full problem using active sets of only a few thousand unknowns. However, it has the drawback that global convergence might fail; see section \ref{S:Counterexample} for a counterexample.
By contrast, for \eqref{eq:ch_many} we can prove global convergence. Unfortunately, this rule has the drawback that it is inefficient in practice due to the huge search space.

Let us now see how much of the  rigorous analysis from section \ref{S:convergence} can be extended to the multi-marginal case. An extension of Theorem \ref{T:cyclelength} on the number of steps needed to find an improving configuration is possible, but a different argument is required because the MMOT equivalent of $c$-cyclical monotonicity fails to provide a practical upper bound on the length of a cycle.

\subsection{Tail clearing for MMOT}
As before, let $\Omega \subset X_1 \times \dots \times X_N$ be a feasible subset of configurations $r=(r_1,...,r_N)$, $\gamma$ the current solution of the reduced problem (i.e. of \eqref{MMOT} with $X_1\times \dots\times X_N$ replaced by $\Omega$), and $u_1,...,u_N$ the current Kantorovich potentials. Assume $\gamma$ is not optimal, then there exists $\tilde \gamma$ with $C[\tilde \gamma] < C[\gamma]$ and $\supp(\tilde \gamma) \leq \sum_i (\ell_i - 1) + 1$, as shown in \cite{friesecke_penka2022}. By construction $\supp(\tilde \gamma) \cap \Omega^c \neq \emptyset$, because otherwise $\gamma$ would not have been a current optimal solution. Hence   there exists $\tilde \Omega \supset \supp(\tilde \gamma)$ with $|\tilde \Omega| < \sum_i \ell_i$.

Define $\tilde \Omega_* = \{r \in \tilde \Omega \colon r \notin \Omega\}$. We claim that there exists $r' \in \tilde \Omega_*$, such that the dual certificate is violated, i.e. $(u_1 \oplus \dots \oplus u_N)(r'):=u_1(r'_1)+\ldots+u_N(r'_N) > c(r')$.

We argue by contradiction. Assume $u_1\oplus\dots\oplus u_N \leq c$ on $\tilde \Omega$. Then,
\begin{align*}
    \int_\Omega c \, d\gamma &= \int_\Omega u_1 \oplus \dots \oplus u_N \, d\gamma = \int_{X_1 \times \dots \times X_N} u_1 \oplus \dots \oplus u_N \, d\gamma \\
    &= \sum_{i = 1}^N \int_{X_i} u_i\, d\mu_i = \int_{X_1 \times \dots \times X_N} u_1 \oplus \dots \oplus u_N \, d\tilde \gamma \\
    &= \int_{\tilde \Omega} u_1 \oplus \dots \oplus u_N \, d\tilde\gamma \leq \int_{\tilde \Omega} c \, d\tilde\gamma.
\end{align*}

Now analogously to the argument for the two marginal case, re-solve the reduced problem on $\Omega \cup \{r'\}$.
Then either 
\[ \inf\{ C[\eta] \colon \eta \in \Pi(\mu_1,\dots,\mu_N),\, \supp(\eta) \subset \Omega \cup \{r'\}\} < C[\gamma], \] or equality holds.
In the latter case, $\gamma$ is still optimal, but the dual solutions $(\tilde u_i)$ must be changed to satisfy the dual certificate $\tilde u_1 \oplus \dots \oplus \tilde u_N \leq c$ on $\Omega \cup \{r'\}$.
By the same argument as before, we again find a configuration in $\tilde \Omega \backslash (\Omega \cup \{r'\})$ violating the dual certificate.
Repeating, after at most $|\tilde \Omega_*|$ steps, all configurations are added and $\tilde \gamma$ is now an  accessible solution, lowering the total cost.

This shows that, also in the multi-marginal case, the size of the reduced problems can be limited by $\beta \cdot (\ell_1+\ldots+\ell_N)$, justifying the tail-clearing procedure and memory efficiency of the algorithm. We summarize this finding in the following theorem.
\begin{thm} \label{T:cyclelength_multi}
    Let $\gamma$ be an extremal optimal solution for the reduced multi-marginal optimal transport problem, and let $(u_i)$ be corresponding Kantorovich potentials.
    Assume $\gamma$ is not optimal for the full MMOT problem.
    Then there exist configurations $r^{(1)}\!,\dots,r^{(k)} \in \Omega^c, k < \sum_i \ell_i$, such that
    \[\inf_{\substack{\tilde \gamma \in \Pi(\mu_1,\dots,\mu_N) \\ \supp(\tilde \gamma) \subset \Omega \cup \{r^{(1)}\!,\dots,r^{(k)}\}}} C[\tilde \gamma] < C[\gamma],\]
    and $r^{(1)}\!,\dots,r^{(k)}$ are accepted by GenCol's acceptance criterion.
\end{thm}
However, the improving configurations $r^{(i)}$ must be proposed by the genetic search rule. Thus we can only obtain a global convergence result for  \eqref{eq:ch_many}.

\begin{cor}
    Suppose $X_1,\dots,X_N$ are discrete spaces of finite cardinality, and the hyperparameter $\beta$ is $\geq 2$. For any marginals, any cost function, and any feasible initial set $\Omega \subset X_1 \times \dots \times X_N$, GenCol with the search rule \eqref{eq:ch_many} converges with probability 1 to an exact solution of the multi-marginal OT problem \eqref{MMOT}.
\end{cor}

\subsection{A Counterexample}\label{S:Counterexample}
Our global convergence proof cannot simply be transferred to the multi-marginal case with the efficient search rule \eqref{eq:ch_single} where children differ from an active configuration by only 1 entry. We present a simple counterexample. Let $X_1=X_2=X_3 = \{1,2,3\}$ and $N=3$. Let further 
\begin{equation*}
  \mu_1 = \mu_2 = \mu_3 := \sum_{x=1}^3\frac{1}{3}\delta_x.
\end{equation*}
The cost function is chosen to be
\begin{align*}
  c(x_1,x_2,x_3) := 
  \begin{cases}
    0 & x_1 = x_2 = x_3\\
    1 & (x_1 \neq x_2) \land (x_1 \neq x_3) \land (x_2 \neq x_3)\\
    2 & \text{else.}
  \end{cases}
\end{align*}
Then the transport plan
\begin{equation*}
  \gamma_0 = \frac{1}{3}(\delta_{(1,2,3)}+\delta_{(2,3,1)}+\delta_{(3,1,2)})
\end{equation*}
is a stationary state for GenCol. Obviously the global optimal plan is
\begin{equation*}
  \gamma^\star = \frac{1}{3}(\delta_{(1,1,1)}+\delta_{(2,2,2)}+\delta_{(3,3,3)}).
\end{equation*}
GenCol proposes new configurations by updating one entry of one active configuration. Independently of the dual solution, all possible configurations that can be  proposed are
\begin{multline*}
    (1,2,2), (1,3,3), (2,2,3), (3,2,3), (2,1,1), (2,2,1), (2,3,2), (2,3,3), \\
    (1,3,1), (3,3,1), (3,2,2), (3,3,2), (3,1,1), (3,1,3), (1,1,2), (2,1,2).
\end{multline*}
The cost for all of them is 2, while the cost for all active configurations in $\gamma_0$ is only 1. Hence for any subset of the configurations listed above added to the active configurations in $\gamma_0$, the optimal plan is again $\gamma_0$. Therefore the current solution $\gamma_0$ never changes and $\gamma^*$ cannot be reached.

The example is designed so that one would have to update two entries of an active configuration to reduce the cost. Any update in just one entry increases the cost. 

Some interesting properties of this example are: 
\begin{itemize} 
    \item The problem is symmetric (i.e., all marginals are equal and the cost is symmetric in its variables), like the Coulomb OT problem arising in electronic structure.
    \item One can replace $\{1,2,3\}$ by a convex independent set (i.e. a set all of whose points are extreme points), in which case the cost $c$ can even be chosen convex.
\end{itemize}

\section{Conclusions}
We rigorously justified the GenCol algorithm in the two-marginal case for arbitrary costs and marginals, showing that it avoids non-minimizing stationary states despite maintaining sparsity. 

For the multi-marginal case, we rigorously justified the algorithm {\it provided} the search rule finds the required configurations described in Theorem \ref{T:cyclelength_multi} with positive probability. Thus GenCol rigorously reduces the storage cost from exponential to linear in the number of marginals. However, the number of search steps might be exponentially large. 

It is an interesting open problem whether the efficient search rule \eqref{eq:ch_single} (which only requires quadratically many search steps in the number of marginals) or any similarly efficient modification can be rigorously justified, 
at least for costs of practical interest like the Coulomb cost or the Wasserstein barycenter cost. 

\bibliographystyle{amsplain}
\bibliography{lit_clean}

\providecommand{\bysame}{\leavevmode\hbox to3em{\hrulefill}\thinspace}
\providecommand{\MR}{\relax\ifhmode\unskip\space\fi MR }
\providecommand{\MRhref}[2]{%
  \href{http://www.ams.org/mathscinet-getitem?mr=#1}{#2}
}
\providecommand{\href}[2]{#2}
\begin{thebibliography}{10}

\bibitem{altschuler2021}
Jason~M. Altschuler and Enric Boix-Adsera, \emph{Hardness results for
  multimarginal optimal transport problems}, Discrete Optim. \textbf{42}
  (2021), 100669.

\bibitem{chandrasekaran2016cutting}
Karthekeyan Chandrasekaran, L{\'a}szl{\'o}~A V{\'e}gh, and Santosh~S Vempala,
  \emph{The cutting plane method is polynomial for perfect matchings},
  Mathematics of Operations Research \textbf{41} (2016), no.~1, 23--48.

\bibitem{li2022maxflow}
Li~Chen, Rasmus Kyng, Yang~P. Liu, Richard Peng, Maximilian~Probst Gutenberg,
  and Sushant Sachdeva, \emph{Maximum flow and minimum-cost flow in
  almost-linear time}, 2022 IEEE 63rd Annual Symposium on Foundations of
  Computer Science (FOCS), 2022, pp.~612--623.

\bibitem{dong2020study}
Yihe Dong, Yu~Gao, Richard Peng, Ilya Razenshteyn, and Saurabh Sawlani, \emph{A
  study of performance of optimal transport}, 2020.

\bibitem{franklin1989scaling}
Joel Franklin and Jens Lorenz, \emph{On the scaling of multidimensional
  matrices}, Linear Algebra and its applications \textbf{114} (1989), 717--735.

\bibitem{friesecke_penka2022}
Gero Friesecke and Maximilian Penka, \emph{The {GenCol} algorithm for
  high-dimensional optimal transport: general formulation and application to
  barycenters and {Wasserstein} splines}, arXiv preprint arXiv:2209.09081
  (2022).

\bibitem{friesecke2022genetic}
Gero Friesecke, Andreas~S Schulz, and Daniela V\"{o}gler, \emph{Genetic column
  generation: Fast computation of high-dimensional multimarginal optimal
  transport problems}, SIAM J. Sci. Comput. \textbf{44} (2022), no.~3,
  A1632--A1654.

\bibitem{lubbecke2005selected}
Marco~E L{\"u}bbecke and Jacques Desrosiers, \emph{Selected topics in column
  generation}, Operations research \textbf{53} (2005), no.~6, 1007--1023.

\bibitem{Santambrogio2015}
Filippo Santambrogio, \emph{Optimal transport for applied mathematicians},
  Birkhäuser, 2015.

\bibitem{scetbon2021low}
Meyer Scetbon, Marco Cuturi, and Gabriel Peyr{\'e}, \emph{Low-rank sinkhorn
  factorization}, International Conference on Machine Learning, PMLR, 2021,
  pp.~9344--9354.

\bibitem{schmitzer2019stabilized}
Bernhard Schmitzer, \emph{Stabilized sparse scaling algorithms for entropy
  regularized transport problems}, SIAM J. Sci. Comput. \textbf{41} (2019),
  no.~3, A1443--A1481.

\bibitem{strossner2022low}
Christoph Str\"{o}ssner and Daniel Kressner, \emph{Low-rank tensor
  approximations for solving multimarginal optimal transport problems}, SIAM J.
  Imaging Sci. \textbf{16} (2023), no.~1, 169--191.

\end{thebibliography}

\end{document}